\documentclass{amsart}

\usepackage{amssymb,amsmath}
\usepackage{mathrsfs}
\usepackage{dsfont}
\usepackage[neveradjust]{paralist}
\usepackage[colorlinks,plainpages,backref]{hyperref}
\usepackage[all]{xypic}
\xyoption{dvips}
\usepackage{url}

\theoremstyle{definition}
\newtheorem{ntn}{Notation}
\newtheorem{dfn}[ntn]{Definition}
\theoremstyle{plain}
\newtheorem{lem}[ntn]{Lemma}
\newtheorem{prp}[ntn]{Proposition}
\newtheorem{thm}[ntn]{Theorem}
\newtheorem{cor}[ntn]{Corollary}

\theoremstyle{remark}

\newcommand{\A}{\mathcal{A}}

\newcommand{\kk}{\mathbf{k}}
\newcommand{\KK}{\mathbb{K}}

\newcommand{\ideal}[1]{{\left\langle#1\right\rangle}}
\newcommand{\into}{\hookrightarrow}
\newcommand{\xymat}{\SelectTips{cm}{}\xymatrix}

\newcommand{\ZZ}{\mathbb{Z}}

\DeclareMathOperator{\Ass}{Ass}
\DeclareMathOperator{\depth}{depth}
\DeclareMathOperator{\Der}{Der}
\DeclareMathOperator{\Hom}{Hom}
\DeclareMathOperator{\pdim}{pdim}
\DeclareMathOperator{\res}{res}
\DeclareMathOperator{\Sym}{Sym}

\begin{document}

\title[Freeness and multirestriction of hyperplane arrangements]
{Freeness and multirestriction of\\hyperplane arrangements}

\author{Mathias Schulze}
\address{
M. Schulze\\
Department of Mathematics\\
Oklahoma State University\\
Stillwater, OK 74078\\
United States}
\email{mschulze@math.okstate.edu}
\thanks{The author gratefully acknowledges support by the ``SQuaREs'' program of American Institute of Mathematics.
He would like to thank Graham Denham, Hal Schenck, Max Wakefield, and Uli Walther for helpful discussions.}

\date{\today}

\begin{abstract}
Generalizing a result of Yoshinaga in dimension $3$, we show that a central hyperplane arrangement in $4$-space is free exactly if its restriction with multiplicities to a fixed hyperplane of the arrangement is free and its reduced characteristic polynomial equals the characteristic polynomial of this restriction.
We show that the same statement holds true in any dimension when imposing certain tameness hypotheses.
\end{abstract}



\maketitle
\tableofcontents

\section{Introduction}

One of the main themes in the study of hyperplane arrangements is the relation of combinatorial and algebraic or topological data associated with arrangements.
In the present article we are concerned with the relation of characteristic polynomials of restrictions on the combinatorial (and topological) side and freeness (of the module of logarithmic differentials/derivations) on the algebraic side.
One of the main open problems in this field, Terao's conjecture, states that freeness is a combinatorial property for simple arrangements.

Let $\A$ be a simple central arrangement of $m$ hyperplanes in an $\ell$-dimensional vector space $V$ over a field $\KK$ of characteristic $0$.
Its intersection lattice $L_\A$, ordered by reverse inclusion, is a geometric lattice with rank function equal to the codimension in $V$.
Note that $V$ is the unique minimal element in $L_\A$.
The M\"obius function of $L_\A$ is the map
\[
\mu\colon L_\A\to\ZZ
\]
defined recursively by the equalities
\begin{align}\label{32}
\nonumber\mu(V)&=1,\\
\forall X\in L_\A\setminus\{V\}\colon\sum_{X\ge Y\in L_\A}\mu(Y)&=0.
\end{align}
It determines the characteristic polynomial of $\A$, defined as 
\[
\chi(\A,t)=\sum_{X\in L_\A}\mu(X)t^{\dim X}.
\]
By \eqref{32}, $\chi(\A,1)=0$ and the quotient
\begin{equation}\label{43}
\chi_0(\A,t)=\chi(\A,t)/(t-1)
\end{equation}
is called the reduced characteristic polynomial.
For complex arrangements, it is equivalent to the Poincar\'e polynomial $\pi(\A,t)=t^\ell\chi(\A,-t^{-1})$ of the complement of the arrangement by \cite{OT80}.

The complex $\Omega^\bullet(\A)$ of logarithmic differential forms along a multiarrangement $\A$ with multiplicities $\kk=(k_H)_{H\in\A}$ consists of the $S$-modules
\[
\Omega^p(\A)=\left\{\omega\in\frac1Q\Omega^p_V\mid\forall H\in\A\colon\frac{d\alpha_H}{\alpha_H^{k_H}}\wedge\omega\in\frac1Q\Omega_V^{p+1}\right\},
\]
and it is closed under exterior product.
The module $\Omega^1(\A)$ and the $S$-module of logarithmic derivations 
\[
D(\A)=D_1(\A)=\{\delta\in D_V\mid\forall H\in\A\colon\delta(\alpha_H)\in\alpha_H^{k_H}S\},
\]
where $D_V=\Der_\KK(S,S)$ denotes the $S$-module of polynomial vector fields on $V$, are mutually $S$-dual.
In particular, $\Omega^1(\A)$ and $D(\A)$ are reflexive $S$-modules.
If one of them is a free $S$-module then $\A$ is called a free arrangement.

Consider now a fixed hyperplane $H\in\A$ and the multiarrangement restriction $\A^H$ of $\A$ to $H$.
If $\A$ is defined by $Q=\prod_{H'\in\A}\alpha_{H'}\in S=\Sym V^*$, and $H$ by $\alpha_H\in V^*$, then $\A^H$ is defined by
\[
Q_H=\frac{Q}{\alpha_H}\vert_H\in S'=\Sym H^*.
\]
For any $H'\in\A\backslash\{H\}$ the natural multiplicity of $H\cap H'\in\A^H$ is 
\[
k_{H\cap H'}=\#\{H''\in\A\backslash\{H\}\mid H\cap H'=H\cap H''\}.
\]

In order to formulate our main result, we extend the well-known notion of tameness for arrangements.

\begin{dfn}
We call $\A$ tame if
\begin{equation}\label{28}
\forall p=1,\dots,\ell-1\colon\pdim\Omega^p(\A)\le p.
\end{equation}
If \eqref{28} holds for $\Omega^p(\A)$ replaced by $D^p(\A)$, we call $\A$ dually tame.
We speak of weak tameness, in both cases, if \eqref{28} holds for $p=1$.
\end{dfn}

Note that all $3$-arrangements are both tame and dually tame.
For locally free arrangements, such as, for example, generic arrangements, weak tameness coincides with ordinary tameness by \cite[Cor.~5.4]{MS01}.

Yoshinaga~\cite[Cor.~3.3]{Yos05} proved our following main result for $\ell=3$.
For the proof we shall use his result \cite[Thm.~2.2]{Yos04} for $\ell>3$ and extend his approach in \cite{Yos05}.

\begin{thm}\label{1}
Assume that $\ell\le4$ or that $\A$ is weakly (dually) tame.
Then $\A$ is free if and only if $\A^H$ is free with exponents $d'_2,\dots,d'_\ell$ and
\begin{equation}\label{2}
\chi_0(\A,t)=\prod_{k=2}^\ell(t-d'_k).
\end{equation}
\end{thm}

The `only if' part holds true for all $\ell\ge3$ by \cite[Thm.~11]{Zie89b} 
and Terao's factorization theorem \cite{Ter81}.

Assuming freeness of $\A^H$, condition~\eqref{2} has the following geometrical interpretation.
Let $\A_s$ be the restriction of $\A$ to the affine space $H_s$ defined by $\alpha_H=s$.
By the natural $\KK^*$-action, the arrangements $\A_s$, $s\ne0$, form a trivial family of non-central simple arrangements and $\A$ can be identified with the cone of $\A_s$, for any $s\ne0$.
Then \cite[Prop.~2.51]{OT92} shows that 
\[
\forall s\ne0\colon\chi(\A_s,t)=\chi_0(\A,t).
\]
On the other hand, the right hand side of \eqref{2} equals the characteristic polynomial $\chi(\A^H)$ of the multiarrangement $\A^H=\A_0$ as defined in \cite{ATW07}.
So condition \eqref{2} can be interpreted as the family $\A_s$ being trivial on the level of characteristic polynomials. 
However, we do not know a concept of characteristic polynomial that covers both the central multiarrangement case and the non-central simple arrangement case.

Whether our result has any implications on Terao's conjecture is, a priori, unclear.
For instance, Ziegler~\cite[Prop.~10]{Zie89b} showed that exponents of multiarrangements are not combinatorial in general.

\section{Multirestriction of logarithmic forms}

Our starting point is a construction from \cite{Yos05}.
Let $M^\bullet$ denote the image of the restriction map
\begin{equation}\label{48}
\res_H^\bullet\colon\Omega^\bullet(\A)\to\Omega^\bullet(\A^H),\quad\frac{d\alpha_H}{\alpha_H}\wedge\eta+\sigma\mapsto\sigma\vert_H,
\end{equation}
and let $C^\bullet$ be its cokernel.
Then $\res_H^\bullet$ factors through $\frac{d\alpha_H}{\alpha_H}\wedge-\colon\Omega^\bullet\to\frac{d\alpha_H}{\alpha_H}\wedge\Omega^\bullet$ into the residue map
\begin{equation}\label{49}
\frac{d\alpha_H}{\alpha_H}\wedge\Omega^\bullet(\A)\to\Omega^\bullet(\A^H).
\end{equation}
As the complex $(\Omega^\bullet(\A),\frac{d\alpha_H}{\alpha_H}\wedge-)$ is exact by \cite[Prop.~4.86]{OT92}, 
\begin{equation}\label{39}
\frac{1}{Q}\Omega^\ell=\Omega^\ell(\A)=\frac{d\alpha_H}{\alpha_H}\wedge\Omega^{\ell-1}(\A)=\frac{d\alpha_H}{\alpha_H}\wedge\frac{\alpha_H}{Q}\Omega^{\ell-1}.
\end{equation}
By definition, 
\begin{equation}\label{40}
\Omega^{\ell-1}(\A^H)=\frac{1}{Q_H}\Omega_H^{\ell-1},
\end{equation}
and it follows using \eqref{39} that $C^{\ell-1}=0$.
Thus, $C^p\ne0$ only for $1\le p\le\ell-2$.

\begin{prp}\label{6}
If $\A^H$ is free and $C^1=0$ then $C^\bullet=0$.
In particular, freeness of $\A$ implies $C^\bullet=0$.
\end{prp}

\begin{proof}
If $\A^H$ is free then, by Saito's criterion \cite[Thm.~8]{Zie89b} and \eqref{40}, 
\[
\Omega^{\ell-1}(\A^H)=\bigwedge^{\ell-1}\Omega^1(\A^H).
\]
So the proof of \cite[Prop.~4.81]{OT92} shows that
\[
\forall p=0,\dots,\ell-1\colon\Omega^p(\A^H)=\bigwedge^p\Omega^1(\A^H).
\]
As restriction commutes with exterior product, these equalities and the surjectivity of $\res_H^1$ implies the surjectivity of $\res_H^\bullet$, and hence $C^\bullet=0$ as claimed.

Now the second statement follows from \cite[Thm.~2.5]{Yos05}
\end{proof}

By \cite[Thm.~2.5, Cor.~2.6]{Yos05} and Terao's factorization theorem, the `if' part of the statement of Theorem~\ref{1} is equivalent to 
\begin{center}
``\eqref{2} implies $C^1=0$ for free $\A^H$''.
\end{center}
The following result is a first step in this direction.

\begin{prp}\label{7}
If $\A^H$ is free and \eqref{2} holds then $\dim(C^p)\le\ell-4$ for all $p\in\ZZ$.
\end{prp}

\begin{proof}
By \cite[Thm.~2.8]{Yos05}, the Poincar\'e series 
\begin{equation}\label{42}
P(M^\bullet,x,y)=\sum_{p=0}^{\ell-1}P(M^p,x)y^p
\end{equation}
of $M^\bullet$ and $\Phi(\A,x,y)$ of $\Omega^\bullet(\A)$ are related by the equality
\begin{equation}\label{33}
\Phi(\A,x,y)=\frac{x+y}{x(1-x)}P(M^\bullet,x,y).
\end{equation}
Assuming freeness of $\A^H$, Ziegler's dual version \cite[Def.-Thm.~3.21]{Zie89a} of \cite{ATW07} and \cite{ST87} shows that 
\begin{align}
\label{34}\chi(\A^H,t)&=\lim_{x\to1}P(\Omega^\bullet(\A^H),x,t(1-x)-1),\\
\label{35}\chi(\A,t)&=\lim_{x\to1}\Phi(\A,x,t(1-x)-1).
\end{align}
As in \cite[(9)]{Yos05}, combining \eqref{43}, \eqref{42}, \eqref{33}, and \eqref{35} yields
\begin{align}\label{36}
\chi_0(\A,t)&=\frac1{t-1}\lim_{x\to1}\frac{x+t(1-x)-1}{x(1-x)}\sum_{p=0}^{\ell-1}P(M^p,x)(t(1-x)-1)^p\\
\nonumber&=\lim_{x\to1}P(M^\bullet,x,t(1-x)-1),
\end{align}
and hence, applying \eqref{34} and \eqref{36} to the exact sequence 
\begin{equation}\label{46}
0\to M^\bullet\to\Omega^\bullet(\A^H)\to C^\bullet\to 0,
\end{equation}
we obtain
\begin{equation}\label{37}
\chi(\A^H,t)-\chi_0(\A,t)=\lim_{x\to1}P(C^\bullet,x,t(1-x)-1).
\end{equation}
Recall that $C^0=0=C^{\ell-1}$ and set $P^p(x)=(-1)^pP(C^p,x)$ for $p=0,\dots,\ell-1$.
Assuming \eqref{2}, the left hand side of \eqref{37} vanishes, and we compute
\begin{align}\label{38}
0&=\lim_{x\to1}P(C^\bullet,x,(-t)(1-x)-1)\\
\nonumber&=\lim_{x\to1}\sum_{p=0}^{\ell-2}(t(1-x)+1)^pP^p(x)\\
\nonumber&=\lim_{x\to1}\sum_{p=0}^{\ell-2}\sum_{r=0}^p{p\choose r}t^r(1-x)^rP^p(x)\\
\nonumber&=\sum_{r=0}^{\ell-2}t^r\lim_{x\to1}\sum_{p=r}^{\ell-2}{p\choose r}(1-x)^rP^p(x).
\end{align}
Considering \eqref{38} as an equality of polynomials in $(\KK[x])[t]$ shows that
\[
\lim_{x\to1}\sum_{p=r}^{\ell-2}{p\choose r}(1-x)^rP^p(x)=0,\quad r=0,\dots,\ell-2.
\]
In particular, for each $k=0,\dots,\ell-2$, we have
\begin{equation}\label{3}
\lim_{x\to1}\sum_{p=r}^{\ell-2}{p\choose r}(1-x)^kP^p(x)=0,\quad r=0,\dots,k.
\end{equation}
For $k=\ell-2$, the matrix of the system of linear equations \eqref{3} has full rank.
By induction on the pole order of $P^\bullet(x)$, it follows that 
\begin{equation}\label{4}
\lim_{x\to1}(1-x)^mP^p(x)=0,\quad p=1,\dots,\ell-2,
\end{equation}
holds for $m=\ell-2$.
Thus, \eqref{3} for $k=\ell-3$ can be considered as a system of linear equations satisfied by $(\lim_{x\to1}(x-1)^{\ell-3}P^p(x))_{p=1,\dots,\ell-2}$ with matrix
\begin{equation}\label{5}
\begin{pmatrix}
{1\choose0} & {2\choose0} & {3\choose0} & \dots & {\ell-3\choose0} & {\ell-2\choose0}\\
{1\choose1} & {2\choose1} & {3\choose1} & \dots & {\ell-3\choose1} & {\ell-2\choose1}\\
0 & {2\choose2} & {3\choose2} & \dots & {\ell-3\choose2} & {\ell-2\choose2}\\
0 & 0 & {3\choose3} & \dots & {\ell-3\choose3} & {\ell-2\choose3}\\
\vdots &  & \ddots & \ddots & \vdots & \vdots\\
0 & \hdots & \hdots & 0 & {\ell-3\choose\ell-3} & {\ell-2\choose\ell-3}
\end{pmatrix}
\end{equation}
Using the relation $\sum_{k=0}^n(-1)^k{n\choose k}=0$, one can eliminate the lower subdiagonal in \eqref{5} by appropriate row operations, keeping all diagonal elements in $\KK^*$.
Thus, \eqref{5} is invertible and \eqref{4} also holds for $m=\ell-3$.
This means that $\ell-4$ is an upper bound for the pole order of $P(C^p,x)$ at $x=1$ which equals $\dim C^p$, for all $p=1,\dots,\ell-2$.
The claim follows.
\end{proof}

The following statement for $\ell=3$ also follows from \cite[Thm.~3.2]{Yos05}.

\begin{cor}\label{12}
Assume that $\A^H$ is free and \eqref{2} holds.
For $\ell=3$, $\dim_\KK C^1=0$ and, for $\ell=4$, $\dim_\KK C^1=\dim_\KK C^2<\infty$.
\end{cor}

\section{Duality of multirestriction for forms and derivations}

To understand $C^{\ell-2}$ and deal with the case $\ell=4$, the following trivial lemma will be useful.
For simple arrangements, it is well known.
We fix a coordinate system $z_1,\dots,z_\ell\in V^*$ and set $dz=dz_1\wedge\dots\wedge z_\ell$.

\begin{lem}\label{14}
Let $\A$ be a multiarrangement in $V$ with multiplicities $\kk=(k_H)_{H\in\A}$ and defining equation $Q\in S$.
\begin{enumerate}[(a)]
\item\label{14a} The inner product induces a $S$-bilinear pairing $D(\A)\times\Omega^p(\A)\to\Omega^{p-1}(\A)$
\item\label{14b} $\ideal{-,dz}$ induces an isomorphism of graded $S$-modules $D(\A)\cong Q\Omega^{\ell-1}(\A)$.
\end{enumerate}
\end{lem}

\begin{proof}\ \pushQED{\qed}
\begin{asparaenum}[(a)]

\item Let $\frac\omega Q\in\Omega^p(\A)$ and $\xi\in D(\A)$.
Then $d\alpha_H\wedge\omega\in\alpha_H^{k_H}\Omega_V^{p+1}$ and $\xi(\alpha_H)\in\alpha_H^{k_H}S$, which implies that
\begin{equation}\label{45}
d\alpha_H\wedge\ideal{\xi,\omega}=\ideal{\xi,d\alpha_H\wedge\omega}+\xi(\alpha_H)\omega\in\alpha_H^{k_H}\Omega_V^p.
\end{equation}
Thus, $\ideal{\xi,\frac{\omega}{Q}}=\frac{\ideal{\xi,\omega}}{Q}\in\Omega^{p-1}(\A)$.

\item By part \eqref{14a}, each $\xi\in D(\A)$ defines an $S$-linear map
\[
\xi=\ideal{\xi,-}\colon S\frac{dz}{Q}=\Omega^{\ell}(\A)\to\Omega^{\ell-1}(\A)
\]
which associates to $\xi$ an element $\frac{\ideal{\xi,dz}}Q=\ideal{\xi,\frac{dz}Q}\in\Omega^{\ell-1}(\A)$.
One computes that 
\begin{equation}\label{44}
\frac{d\alpha_H}{\alpha_H^{k_H}}\wedge\ideal{\xi,\frac{dz}{Q}}=\frac{\xi(\alpha_H)}{\alpha_H^{k_H}}dz
\end{equation}
which shows that $\xi\mapsto\ideal{\xi,\frac{dz}{Q}}$ is an isomorphism.\qedhere

\end{asparaenum}
\end{proof}

For the simple arrangement $\A$, consider the module of derivations that are logarithmic both along $\A$ and along the level sets of $\alpha_H$,
\[
D_H(\A)=D(\A)\cap D(\alpha_H)=\{\delta\in D(\A)\mid\delta(\alpha_H)=0\}.
\]
Then, by \cite[Thm.~11]{Zie89b}, there is a restriction map
\[
\res_H\colon D_H(\A)\to D(\A^H),
\]
and we denote by $M$ its image and by $C$ its cokernel.

\begin{prp}\label{13}
There is an isomorphism $C[m-1]\cong C^{\ell-2}$ of graded $S'$-modules.
\end{prp}

\begin{proof}
By Lemma~\ref{14}.\eqref{14b}, there is an isomorphism of graded $S$-modules 
\begin{equation}\label{9}
\ideal{-,dz}\colon D(\A)\to Q\Omega^{\ell-1}(\A).
\end{equation}
Using the Euler derivation $\chi\in D(\A)$, one can decompose
\begin{equation}\label{15}
D(\A)=S\chi\oplus D_H(\A).
\end{equation}
Then \eqref{44} shows that $\ideal{D_H(\A),\frac{dz}{Q}}\subset\Omega^{\ell-1}(\A)$ is the kernel of $\frac{d\alpha_H}{\alpha_H}\wedge-$.
Using that the complex $(\Omega^\bullet(\A),\frac{d\alpha_H}{\alpha_H}\wedge-)$ is exact by \cite[Prop.~4.86]{OT92}, we conclude that
\[
\ideal{D_H(\A),\frac{dz}{Q}}=\frac{d\alpha_H}{\alpha_H}\wedge\Omega^{\ell-2}(\A).
\]
As $M^{\ell-2}$ is the residue of the latter module along $H$, \eqref{9} induces an isomorphism 
\begin{equation}\label{10}
M=D_H(\A)\vert_H\cong Q_HM^{\ell-2}.
\end{equation}
Again by Lemma~\ref{14}.\eqref{14b}, there is an isomorphism of $S'$-modules
\begin{equation}\label{11}
D(\A^H)\cong Q_H\Omega^{\ell-2}(\A^H).
\end{equation}
The claim follows by combining \eqref{10} and \eqref{11}.
\end{proof}

Assume that $\ell=4$, that $\A^H$ is free, and that \eqref{2} holds.
By Corollary~\ref{12}, Proposition~\ref{13}, and homogeneity, the cokernel $C$ of $\res_H$ is supported at $0\in H$ only.
Therefore, $\A$ is locally free along $H$ by \cite[Thm.~2.1]{Yos04}.
Then \cite[Thm.~2.2]{Yos04} yields the following result.

\begin{thm}\label{41}
The statement of Theorem~\ref{1} holds true for $\ell=4$.
\end{thm}

\section{Tameness hypotheses on forms and derivations}

In order to apply the (weak) tameness hypothesis to \eqref{48} and \eqref{49}, we need the following trivial result.

\begin{lem}\label{29}
There is a direct sum decomposition 
$\Omega^1(\A)\cong S\frac{d\alpha_H}{\alpha_H}\oplus\frac{d\alpha_H}{\alpha_H}\wedge\Omega^1(\A)$.
\end{lem}

\begin{proof}
Equality \eqref{45} with $\xi=\chi$ shows that 
\[
0\to S\frac{d\alpha_H}{\alpha_H}\to\Omega^1(\A)\to\frac{d\alpha_H}{\alpha_H}\wedge\Omega^1(\A)\to0
\]
is split exact, cf.\ the proof of \cite[Prop.~4.86]{OT92}.
\end{proof}

To see the relation with the direct sum decomposition \eqref{15}, consider the split exact sequence
\begin{equation}\label{31}
\xymat{
0\ar[r] & D_H(\A)\ar[r] & D(\A)\ar[r]^-\varphi & S\ar[r] & 0
}
\end{equation}
where
\begin{equation}\label{30}
\varphi(\delta)=\frac{\delta(\alpha_H)}{\alpha_H}=\ideal{\frac{d\alpha_H}{\alpha_H},\delta}.
\end{equation}
Applying $-^\vee=\Hom_S(-,S)$ to \eqref{31} yields a split exact sequence
\begin{equation}\label{47}
\xymat{
0 & D_H(\A)^\vee\ar[l] & \Omega^1(\A)\ar[l] & S\ar[l]_-{\varphi^\vee} & 0\ar[l].
}
\end{equation}
By \eqref{30} and \eqref{47}, we can identify $\varphi^\vee=\frac{d\alpha_H}{\alpha_H}$, and hence
\[
D_H(\A)^\vee\cong\Omega^1(\A)/S\frac{d\alpha_H}{\alpha_H}\cong\frac{d\alpha_H}{\alpha_H}\wedge\Omega^1(\A).
\]

\begin{thm}\label{27}
The statement of Theorem~\ref{1} holds true if $\A$ is weakly (dually) tame.
\end{thm}

\begin{proof}
By \cite[Cor.~3.3]{Yos05} and Theorem~\ref{41}, we may assume that $\ell\ge5$.
We consider the weakly tame case only; the weakly dually tame case can be treated along the same lines using the direct sum decomposition \eqref{15} instead of Lemma~\ref{29}, and proving surjectivity of $\res_H$ instead of $\res_H^1$.

By Lemma~\ref{29} and weak tameness, $\pdim_S(\frac{d\alpha_H}{\alpha_H}\wedge\Omega^1(\A))\le1$.
Then also $\pdim_{S'}M^1\le1$ due to the exact sequence
\[
\xymat{
0\ar[r] & \frac{d\alpha_H}{\alpha_H}\wedge\Omega^1(\A)\ar[r]^-{\alpha_H} & \frac{d\alpha_H}{\alpha_H}\wedge\Omega^1(\A)\ar[r] & M^1\ar[r] & 0.
}
\]
Indeed, if $F_\bullet\colon F_1\into F_0$ is an $S$-free resolution of $\frac{d\alpha_H}{\alpha_H}\wedge\Omega^1(\A)$, then $F_\bullet/\alpha_H F_\bullet$ is an $S'$-free resolution of $M^1$.
This follows immediately from comparing the two spectral sequences associated to the double complex defined by multiplication by $\alpha_H$ on $F_\bullet$.

By the Auslander--Buchsbaum theorem, this shows that
\begin{equation}\label{16}
\forall x\in H\colon\dim S'_x\ge3\Rightarrow\depth M^1_x\ge 2,
\end{equation}
where $x$ is considered as a scheme-theoretic point.
From \eqref{46} we get an exact sequence of $S'$-modules
\begin{equation}\label{25}
0\to M^1\to\Omega^1(\A^H)\to C^1\to0.
\end{equation}
Assuming that $\A^H$ is free, \eqref{16} applied to the localization of \eqref{25} at $x$ shows that $\depth C^1_x\ge1$ for any $x\in H$ of codimension at least $3$.
However, by Proposition~\ref{7} and our hypothesis $\ell\ge5$, $C^1$ has codimension at least $3$ in $H$.
Thus, $\Ass C^1=\emptyset$, $C^1=0$, $\res^1_H$ is surjective, and $\A$ is free by \cite[Thm.~2.5]{Yos05}.
\end{proof}

\bibliographystyle{amsalpha}
\bibliography{ycor}
\end{document}